\documentclass{article}

\newtheorem{theorem}{Theorem}
\newtheorem{definition}{Definition}
\newtheorem{example}{Example}
\newtheorem{corollary}{Corollary}

\def\HH{ { H} }
\def\pd#1{ \partial_{#1} }
\title{Modified ${\cal A}$-hypergeometric Systems}
\author{
Nobuki Takayama,
Department of Mathematics, Kobe University
}
\date{June 30, 2007, January 20, 2008 revised}

\begin{document}
\maketitle

\noindent
{\bf Abstract}
We will introduce a modified system of ${\cal A}$-hypergeometric system
(GKZ system) by applying a change of variables for Gr\"obner deformations
and study its Gr\"obner basis and the indicial polynomial along the 
exceptional hypersurface.

\section{Introduction}

Since the work of Gel'fand, Zelevinsky, and Kapranov \cite{GZK},
studies of ${\cal A}$-hypergeometric system (GKZ system) have 
attracted a lot of mathematicians, who want to understand
hypergeometric differential equations in a general way.
We refer to the book \cite{SST} on the status of the art in 2000,
and the recent papers \cite{mmw} and \cite{SW} and their reference trees
on recent advances.
We also note that these studies have had  fruitful interactions
with frontiers of computational commutative algebra and
computational $D$-modules.

In this short paper, we will introduce a modified version of
this ${\cal A}$-hypergeometric system
and provide a first step to study it.
The original system is defined on the $y=(y_1, \ldots, y_n)$ space
and the modified system is defined on the $(t,x_1,\ldots, x_n)$ space
with one more variable $t$.
Let us sketch our idea to introduce the modified system.
We consider the direct sum of the ${\cal A}$-hypergeometric system on the $y$ space
and the $D$-module $D/D\cdot s \pd{s}$ on the $s$-space.
For a weight vector $w \in {\bf Z}^n$, 
the original system restricted on the complex torus
is transformed into the modified system on $(t,x)$ space by the map
$$ {\bf C}^n \times {\bf C}^* \ni (y_1, \ldots, y_n,s) \mapsto
   (t^{w_1} x_1, \ldots, t^{w_n} x_n,t) \in {\bf C}^n \times {\bf C}^* $$
(see \cite{Sturmfels-kyushu} and \cite{SST} on this transformation).
The transformed system can be naturally extended on ${\bf C}^{n+1}$.
Intuitively speaking, the variety $t=0$ is analogous to the exceptional hypersurface
of a blowing-up operation.
We will study the indicial polynomial along $t=0$ as a first step
to make a local and global analysis of the modified system.
As a byproduct of our discussion on the modified system,
we will also give a proof to the claim
${\rm rank}(H_A(\beta)) \geq {\rm vol}(A)$ for non-homogeneous $A$.

\section{Definition and Holonomic Rank of Modified ${\cal A}$-hypergeometric systems}
Let $A=(a_{ij})_{ij}$ be a $d\times n$-matrix whose elements are integers
and $w=(w_1, \ldots, w_n)$ a vector of integers. 
We suppose that the set of the column vectors of $A$ spans ${\bf Z}^d$.
Define 
$${\tilde A} = 
  \pmatrix{ a_{11} & \cdots & a_{1n} & 0 \cr
                   & \cdots &        & 0 \cr
            a_{dn} & \cdots & a_{dn} & 0 \cr
            w_1    & \cdots & w_n    & 1 \cr}.
$$                 

\begin{definition} \rm
We call the following system of differential equations $\HH_{A,w}(\beta)$
a {\it modified  ${\cal A}$-hypergeometric differential system}:
\begin{eqnarray*}
  \left( \sum_{j=1}^n a_{ij}  x_j \partial_{j} - \beta_i \right) \bullet f &=& 0,
   \qquad(i = 1, \ldots, d)  \\
  \left( \sum_{j=1}^n w_j x_j \partial_{j} - t \partial_t \right) \bullet f &=& 0,
    \\
  \left( \prod_{i=1}^n \partial_{i}^{u_{i}} t^{u_{n+1}}
       - \prod_{j=1}^n \partial_{j}^{v_{j}} t^{v_{n+1}}
  \right) \bullet f &=& 0.
   \qquad( u, v \in {\bf N}^{n+1} \mbox{ and } {\tilde A}u = {\tilde A}v )
\end{eqnarray*}
\end{definition}

\medbreak
\noindent
Let $I_{\tilde A}$ be the toric ideal 
generated by
\begin{equation}
   \prod_{i=1}^n \partial_{i}^{u_{i}} t^{u_{n+1}}
       - \prod_{j=1}^n \partial_{j}^{v_{j}} t^{v_{n+1}}
   \qquad( u, v \in {\bf N}^{n+1} \mbox{ and } {\tilde A}u = {\tilde A}v )
\end{equation}
in ${\bf C}[\pd{1}, \ldots, \pd{n}, t]$.
Since
${\bf C}[\pd{1}, \ldots, \pd{n}, t]/I_{\tilde A}$ is an integral domain
and $t^m$ does not belong to the toric ideal,
we have
\begin{equation}
 I_{\tilde A} = I^{sat}_{\tilde A} = (I_{\tilde A} : t^\infty)
                    = \{ \ell \,|\, t^m \ell  \in I_{\tilde A} \ 
 \mbox { for a non-negative integer $m$} \} 
\end{equation}
This fact will be used in the proof of Theorem \ref{th:gb}.

We note that the matrix $\tilde A$ with $w=(1, \ldots, 1)$
was introduced in \cite{ot} to construct ${\rm vol}(A)$ 
convergent series solutions.

Throughout this paper, we will use notations and facts shown in \cite{SST}.
In particular, we do not cite original papers
for text level well-known facts in the theory of D-modules.
Refer references of \cite{SST} as to these original papers.

Let $a_i$ be the $i$-th column vector of the matrix $A$ and $F(\beta,x,t)$
the integral
\[
F(\beta,x,t) = \int_C \exp\left(\sum_{i=1}^n x_i t^{w_i} s^{a_i} \right)
s^{-\beta-1} ds, \qquad s=(s_1, \ldots, s_d), \ 
\beta = (\beta_1, \ldots, \beta_d).
\]

The integral $F(\beta,x,t)$ satisfies the modified ${\cal A}$-hypergeometric
differential system ``formally''.
We use the word ``formally''
because, there is no general and rigorous description about the cycle $C$.
However, the integral representation gives an intuitive figure of
what are solutions of modified ${\cal A}$-hypergeometric systems.
The proof is analogous to \cite[221--222]{SST}.
We note that 
if $a_{di}=1$ for all $i$, 
we also have the following ``formal'' integral representation
\begin{eqnarray*}
F(\beta,x,t) &=& \int_C \left(\sum_{i=1}^n x_i t^{w_i} {\tilde s}^{{\tilde a}_i} \right)^{-\beta_d} 
{\tilde s}^{-{\tilde \beta}-1} d{\tilde s}, \\
& & {\tilde a}_i = (a_{1i}, \ldots, a_{d-1,i})^T, 
{\tilde s}=(s_1, \ldots, s_{d-1}), \ 
\beta = (\beta_1, \ldots, \beta_d).
\end{eqnarray*}

We denote by $D$ 
the ring of differential operators
${\bf C}\langle
x_1, \ldots, x_n, t,
\partial_1, \ldots, \partial_n, \partial_t \rangle$.
We will regard modified ${\cal A}$-hypergeometric system
as the left ideal in $D$.
We will denote by $H_{A,w}(\beta)$ the left ideal
as long as no confusion arises.

\begin{theorem}\label{thm:1}
\begin{enumerate}
\item The left $D$-module $D/\HH_{A,w}(\beta)$ is holonomic.
\item The rank of $\HH_{A,w}(\beta)$ agrees with the holonomic rank of $\HH_A(\beta)$
for any $w$.
\end{enumerate}
\end{theorem}

{\it Proof}.
(1) We apply the {\it Laplace transformation} with respect to the variable
$t$ ( $t \mapsto   - \pd{t'}$, $ \pd{t} \mapsto t'$ )
for the modified ${\cal A}$-hypergeometric system $H_{A,w}(\beta)$.
Then, the transformed system is nothing but ${\cal A}$-hypergeometric system
for the matrix ${\tilde A}$ and
the parameter vector $(\beta_1, \ldots, \beta_n, -1)$.
It is known that the transformed system is holonomic, then
the original system is also holonomic by showing the Hilbert polynomials 
with respect to the Bernstein filtration of each system
agree.

\bigbreak

\noindent
(2)
We consider the biholomorphic map  $\varphi$
on ${\bf  C}^n \times {\bf C}^*$
\begin{equation}
{\bf  C}^n \times {\bf C}^* \ni (y_1, \ldots, y_n,s) \mapsto
   (t^{w_1} x_1, \ldots, t^{w_n} x_n,t) \in {\bf  C}^n \times {\bf C}^*
\end{equation}
The map $\varphi$ induces a correspondence of differential operators
on ${\bf  C}^n \times {\bf C}^*$
\begin{eqnarray*}
  \frac{\partial}{\partial y_i} &=& t^{-w_i} \frac{\partial}{\partial x_i} \\
  -s \frac{\partial}{\partial s} &=& 
   -t \frac{\partial}{\partial t} + \sum_{j=1}^n w_n x_n \frac{\partial}{\partial x_i}
\end{eqnarray*}
Consider a left ideal $H_Y$ in 
$D_Y={\bf C}\langle y_1, \ldots, y_n, s, \pd{y_1}, \ldots, \pd{y_n}, \pd{s} \rangle$
generated by $H_A(\beta)$ and  $ s \pd{s}$.
The holonomic rank of $D_Y/H_Y$ is that of $H_A(\beta)$.
We can see that 
the image of ${\cal D}_Y/{\cal D}_Y H_Y$ by the biholomorphic map $\varphi$
on ${\bf C}^n \times {\bf C}^*$ is ${\cal D}_X/{\cal D}_X H_{A,w}(\beta)$
by utilizing the correspondence of differential operators.
Here, ${\cal D}_X$ and ${\cal D}_Y$ denote the sheaves of differential operators
on ${\bf C}^n \times {\bf C}^*$ of $(y,s)$-space and $(x,t)$-space respectively.
Since the holonomic rank agrees with the multiplicity of the zero section
of the characteristic variety at generic points, 
the holonomic ranks of the both systems agree
\cite[pp 28--40]{SST}.  Q.E.D.
\bigbreak

\begin{corollary}  \label{cor:rank}
${\rm rank}\,(H_A(\beta)) \geq {\rm vol}(A)$
\end{corollary}

{\it Proof}.
When $A$ has $(1,1,\ldots, 1)$ in its row space ($A$ is homogeneous),
${\rm rank}\,(H_A(\beta)) \geq {\rm vol}(A)$ holds
\cite[Theorem 3.5.1]{SST}, which is proved by utilizing that $H_A(\beta)$
is regular holonomic and by constructing ${\rm vol}(A)$ many series solutions.
Put $w=(1,1,\ldots,1)$ in the modified system $\HH_{\tilde A}(\beta)$.
Then, we have ${\rm rank}\,(\HH_{\tilde A}(\beta)) \geq {\rm vol}(A)$.
Hence, Theorem \ref{thm:1} gives the conclusion.
Q.E.D.

Note that the upper semi continuity theorem of holonomic rank of
\cite{mmw} also gives this result.

\begin{example} \rm \label{example:1}
We take $A=(1,3)$, $\beta=(-1)$, and $w=(-1,0)$
({\it Airy type integral}) \cite[p.223]{SST}.
Define a sequence $d_m$ by
$$ d_0 = 1, d_{m+1} = \frac{- (3m+1)(3m+2)(3m+3)}{m} d_m $$
The divergent series
\begin{eqnarray}
  f(x;t)&=&\sum_{m=0}^\infty \left( d_m x_1^{-3m-1} x_2^m\right) t^{3m+1} \nonumber \\
       &=& \sum_{m=0}^\infty \left( \frac{\Gamma(3m+1)}{\Gamma(m+1)}
                             x_1^{-3m-1} x_2^m  \right) t^{3m+1} \label{eq:airysol} 
\end{eqnarray}
is a formal solution of the modified system.
Fix a point $(x_1,x_2)=(a_1,a_2)$ such that $a_1,a_2 \not= 0$.
Then this is a Gevrey formal power series
solution at $(a_1,a_2,0)$ along $t=0$ in the class $s=1+2/3$
from the definition of Gevrey series.
The slope of this system can be computed by our program \cite[command {\tt sm1.slope}, {\tt slope}]{openxm}, \cite{ACG},
\cite{ct1} 
and the set of the slopes is $\{-3/2\}$.  
Since $1/(1-s)$ is the slope, 
we have constructed a formal power series standing for the slope.
\end{example}

\section{A Gr\"obner Basis of  Modified ${\cal A}$-Hypergeometric Systems}

We will call $t=0$ the {\it exceptional hypersurface} and we are interested in
local analysis near $t=0$.
We denote by $\tau=({\bf 0},-1; {\bf 0},1)$
the weight vector such that $t$  has the weight $-1$ and
$\pd{t}$ has the weight $1$.
We also denote by ${\tilde A}_{\theta,w,\beta}$
the first $(d+1)$ Euler operators of the modified ${\cal A}$-hypergeometric
system.

It is easy to see that,
for generic $w$, ${\rm in}_\tau(D \cdot I_{\tilde A})$ is generated
by monomials in ${\bf C}[\pd{1}, \ldots, \pd{n}]$
and we will regard it as a monomial ideal in  this commutative ring.

\begin{theorem} \label{th:gb}
For generic $\beta$ and $w$, we have
\begin{equation}
 {\rm in}_{({\bf 0},-1; {\bf 0},1)}(H_{A,w}(\beta)) =
  D \cdot {\rm in}_\tau(D \cdot I_{\tilde A}) 
+ D \cdot {\tilde A}_{\theta,w,\beta}
\end{equation}
\end{theorem}

{\it Proof}.
The proof is analogous to \cite[Theorem 3.1.3]{SST}.
Let $s=(s_1, \ldots, s_d)$ be a vector of new indeterminates.
Consider the algebra 
$$D[s] = {\bf C}\langle x_1, \ldots, x_n, t, \pd{1}, \ldots, \pd{n},\pd{t}, s_1, \ldots, s_d \rangle$$
and its homogenized Weyl algebra by $h$ $D[s]^h$.
Let $H$ be the left ideal in $D[s]^h$ generated by ${\tilde A}_{\theta,w,s^2}$ and
the homogenization of $I_{\tilde A}$.
We define a partial order $>_\tau$ on monomials in $D[s]$ by
\begin{eqnarray*}
 s^a x^b \pd{}^c t^d \pd{t}^{e} >_\tau s^{a'} x^{b'} \pd{}^{c'} t^{d'} \pd{t}^{e'}
&\Leftrightarrow &  -d+e > -d'+e', \ \mbox{or} \\ 
& & -d+e = -d'+e' \ \mbox{and}\ (a,e,d) >_{lex} (a',e',d')
\end{eqnarray*}
We refine this partial order by any monomial order and define orders $<$ in $D[s]$.
(This order on $D[s]$ is extended to the order in the homogenized Weyl algebra 
 and $D[s]^h$ as in \cite[Chapter 1]{SST}.)

Let ${\cal G}$ be the reduced Gr\"obner basis of the homogenized binomial ideal $I_{\tilde A}$ in $D[s]^h$
with respect to the order $<$.
Note that the reduced Gr\"obner basis consists of elements of the form
$\underline{\pd{}^u h^{p}}-\pd{}^v t^{v_{n+1}} h^{p'}$, $v_{n+1}>0$
because $w$ is generic and $I_{\tilde A}$ is saturated with respect to $t$.
Note that either $p=0$ or $p'=0$ holds.

We will show that ${\cal G}$ and ${\tilde A}_{\theta,w,s^2}^h$ is a Grobner basis ${\cal G}'$ with respect to $<$
in $D[s]^h$.
This fact can be shown by checking the S-pair criterion in $D[s]^h$.
It is easy to see that
$$ sp(\underline{\theta_t}- \sum w_j \theta_j, \underline{s_i^2} - \sum a_{ij} \theta_j)
  \rightarrow_{{\cal G}'} 0 $$
$$ sp(\underline{s_k^2} - \sum a_{kj} \theta_j,\underline{s_i^2} - \sum a_{ij} \theta_j)
  \rightarrow_{{\cal G}'} 0 $$
We assume $p>0$ and $p'=0$.
\begin{eqnarray*}
& & sp(\underline{\pd{}^u h^p} - \pd{}^v t^{v_{n+1}}, \underline{s_i^2} - \sum a_{ij} \theta_j) \\
&=& s_i^2(\underline{\pd{}^u h^p} - \pd{}^v t^{v_{n+1}})-
   \pd{}^u h^{p} (\underline{s_i^2} - \sum a_{ij} \theta_j)) \\
&=&-s_i ^2 \pd{}^v t^{v_{n+1}} + \pd{}^u h^{p} \sum a_{ij} \theta_j \\
&=&-s_i ^2 \pd{}^v t^{v_{n+1}} + \left( \sum a_{ij} \theta_j\right) \pd{}^u h^{p} 
                             + \left( \sum a_{ij} u_j \right) \pd{}^u h^{p}  \\
& & \mbox{since $\pd{}^u h^{p} > \pd{}^v t^{v_{n+1}}$ we may rewrite it as} \\
&=& -s_i^2 \pd{}^v t^{v_{n+1}} + \left( \sum a_{ij} \theta_j\right) (\pd{}^u h^{p}- \pd{}^v t^{v_{n+1}})
       + \left(\sum a_{ij} \theta_j \right) \pd{}^v t^{v_{n+1}} + \left(\sum a_{ij} u_j\right)\pd{}^u h^p  \\
&=&  \left( \sum a_{ij} \theta_j \right) (\pd{}^u h^{p}- \pd{}^v t^{v_{n+1}})
    + \pd{}^v t^{v_{n+1}} \left( \sum a_{ij} \theta_j - \sum a_{ij} v_j - s_i^2 \right) 
    + \left( \sum a_{ij} u_j \right) \pd{}^u h^p  \\
& & \mbox{ since $\sum a_{ij} u_j = \sum a_{ij} v_j$ } \\
&=&  \left( \sum a_{ij} \theta_j \right) (\underline{\pd{}^u h^{p}}- \pd{}^v t^{v_{n+1}})
       + \pd{}^v t^{v_{n+1}} \left(\sum a_{ij} \theta_j- s_i^2\right) +
        \left( \sum a_{ij} u_j \right) ( \pd{}^u h^p - \pd{}^v t^{v_{n+1}} ) \\
&\rightarrow_{{\cal G}'}& 0 
\end{eqnarray*}
The case $p=0$, and $p'>0$ can be shown analogously.

The final case we have to check is that
\begin{eqnarray*}
& & sp(\underline{\pd{}^u h^p} - \pd{}^v t^{v_{n+1}}, \underline{\theta_t} - \sum a_{ij} \theta_j) \\
&=& -\theta_t \pd{}^v t^{v_{n+1}} + h^p \pd{}^u \sum w_j \theta_j \\
&=& -\theta_t \pd{}^v t^{v_{n+1}} + \left( \sum w_j \theta_j + \sum w_j u_j\right) h^p \pd{}^u \\
&=& -\theta_t \pd{}^v t^{v_{n+1}} + \left( \sum w_j \theta_j + \sum w_j u_j\right) 
                                   \left(h^p \pd{}^u - \pd{}^v t^{v_{n+1}} \right) \\
& & \quad\quad\quad + \left( \sum w_j \theta_j + \sum w_j u_j\right) \pd{}^v t^{v_{n+1}} \\
&=& -\theta_t \pd{}^v t^{v_{n+1}} + \left( \sum w_j \theta_j + \sum w_j u_j\right) 
                                   \left(h^p \pd{}^u - \pd{}^v t^{v_{n+1}} \right) \\
& & \quad\quad +\pd{}^v t^{v_{n+1}} \left( \sum w_j \theta_j + \sum w_j u_j - \sum w_j v_j \right)  \\
&=&  \left( \sum w_j \theta_j + \sum w_j u_j\right) 
     \left(h^p \pd{}^u - \pd{}^v t^{v_{n+1}} \right) \\
& & \quad +\pd{}^v t^{v_{n+1}} \left( \sum w_j \theta_j + \sum w_j u_j - \sum w_j v_j 
                                     -\theta_t - v_{n+1}\right)  \\
&=&  \left( \sum w_j \theta_j + \sum w_j u_j\right) 
     \left(\underline{h^p \pd{}^u} - \pd{}^v t^{v_{n+1}} \right) 
    +\pd{}^v t^{v_{n+1}} \left( \sum w_j \theta_j -\underline{\theta_t} \right)  \\
&\rightarrow_{{\cal G}'}& 0 
\end{eqnarray*}
The rests of the proof are analogous to \cite[p.106]{SST}.
Q.E.D.

\section{Indicial Polynomial along $t=0$}

We fix generic $w$.
Let $M$ be the monomial ideal ${\rm in}_\tau(I_{\tilde A})$
in ${\bf C}[\pd{1}, \ldots, \pd{n}]$.
The  top dimensional standard pairs are denoted by
${\cal T}(M)$ \cite[p.112]{SST}
and $\beta^{(\pd{}^\beta,\sigma)}$ is the zero point in ${\bf C}^n$
of the distraction of $M$ and $A \theta-\beta$ associated 
to the standard pair $(\pd{}^\beta,\sigma)$.

\begin{theorem}  \label{thm:indicial}
Let $\beta$ and $w$ both be generic.
Then, the indicial polynomial ($b$-function) of $H_{A,w}(\beta)$ 
along $t=0$ is
\begin{equation} \label{eq:indicial}
\sum_{(\pd{}^\beta,\sigma) \in {\cal T}(M)} (s - w \cdot \beta^{(\pd{}^\beta,\sigma)})
\end{equation}
If ${\cal T}(M)$ is the empty set, the indicial polynomial is $0$.
\end{theorem}

{\it Proof}.
Under Theorem \ref{th:gb},
the proof is analogous to \cite[p.198, Proposition 5.1.9]{SST}.
\bigbreak

If the indicial polynomial is not zero and the difference of roots are not integral,
we can construct formal series solution of the form
\begin{equation}  \label{eq:gseries}
 t^e \sum_{k=0}^\infty c_k(x) t^k, \quad c_k \in {\bf C}[1/x_1, \ldots, 1/x_n, x_1, \ldots, x_n] 
\end{equation}
where $e$ is a root of the indicial polynomial and
$t^e c_0(x)$ is a solution of the initial system ${\rm in}_{({\bf 0},-1; {\bf 0},1)} ( H_{A,w}(\beta) )$.
If the indicial polynomial is zero, there is no formal series solution of the form above.

\begin{example} \rm
(Continuation of Example \ref{example:1}).
Note that ${\rm in}_\tau(I_{\tilde A}) = \langle \pd{2} \rangle$.
The distraction \cite[p.68]{SST} 
of ${\rm in}_\tau(H_{A,w}(\beta))$
is generated by
$\theta_2, \theta_1+3 \theta_2+1 , -\theta_1-\theta_t$.
Therefore, the set of zero points are
$\{ (-1,0,1) \}$.
Then, the indicial polynomial is $s-1$.
The formal solution (\ref{eq:airysol}) stands for the root $s=1$.
\end{example}

\begin{example} \rm
Consider the modified hypergeometric system for
$A=(-1,1,2)$, $\beta=(1/2)$, $w=(-2,-1,0)$.
This is the {\it Bessel function in two variables} called by Kimura and Okamoto \cite{ost}.
Although it is a side story in view of this paper, we want to note that
a 3-D Graph of a solution of this system can be seen
at \\
{\tt http://www.math.kobe-u.ac.jp/HOME/taka/test-bess2m.html}.
You will be able to see waves in two directions. 

The indicial polynomial is $0$, because
$I_{\tilde A} \ni \underline{1}-\pd{1}^2 \pd{3}$.
Then, there exists no series solution of the form (\ref{eq:gseries}).  
Incidentally, the set of the slopes along $t=0$ at $x=(2,2,1)$ is equal to
$\{ -2, -3/2 \}$. The values are obtained by our program \cite{openxm}.

Let us change $w$ into $w=(3,2,1)$.
The set of zero points of the distraction are \\
$\{ (-1/2,0,0), (0,0,0), (0,1,-1/4)\}$
(which are obtained by computing 
the primary decomposition of the ideal generated by
the distraction of ${\rm in}_\tau(I_{\tilde A})$ and 
${\tilde A}_{\theta,w,\beta}$)
and the indicial polynomial is $(s-3/2)(s+1/4)(s-7/4)$
(we use the {\tt Risa/Asir} command \verb@ generic_bfct @).
In this case, the generic condition for the Theorem \ref{thm:indicial}
satisfied and the equality (\ref{eq:indicial}) also gives the same answer.
Incidentally,
the local monodromy group of the local solutions around $t=0$ 
is generated by
${\rm diag}(-1, \exp(\pi \sqrt{-1}/2), -\exp(\pi \sqrt{-1}/2))$.
The set of the slopes along $t=0$ at $x=(2,2,1)$ is empty.
\end{example}

It is an interesting open problem to construct rank many series solutions
in terms of formal puiseux series and exponential functions along
$t=0$. 

\bigbreak

{\it Acknowledgements}:
This work was motivated by comments by Bernd Sturmfels
to the joint work with Francisco
Castro-Jimenez on slopes for ${\cal A}$-hypergeometric system \cite{ct1}.
He said that ``can you  construct series solutions standing for the slopes''?
We realized that the original ${\cal A}$-hypergeometric system
has few classical solutions standing for slopes for some examples.
The author introduced  a modified system, which seems to be easier to analyze 
than the original system
and may be a step to study the original ${\cal A}$-hypergeometric system,
when Castro stayed in Japan in the spring of 2006.
In fact, our Example \ref{example:1} gives an example to his question.
However, we are far from a complete answer.
The author thanks to all comments and discussions with them.

The author also thanks to Christine Berkesch
who made me some questions on the first versioin of this paper, 
which yields a substantial improvement in the second version.

Finally, Go Okuyama made a question on the lower bound of the rank
of ${\cal A}$-hypergeometric system. 
The Corollary \ref{cor:rank} is an answer to his question.

\bibliographystyle{plain}

\end{document}